\documentclass[12pt]{amsart}
\usepackage {amssymb,amsmath}

\makeatletter
\def\@cite#1#2{{\m@th\upshape\bfseries%
[{#1\if@tempswa{\m@th\upshape\mdseries, #2}\fi}]}} \makeatother
\theoremstyle{plain}
\newtheorem{thm}{Theorem}[section]

\newtheorem{prop}[thm]{Proposition}

\theoremstyle{definition}
\newtheorem{rem}[thm]{Remark}

\newtheorem{defn}[thm]{Definition}
\newtheorem{note}[thm]{Note}

\newtheorem{eg}[thm]{Example}
\newtheorem{egs}[thm]{Examples}

\newcommand{\Prf}{\noindent\textbf{Proof.\ }}
\newcommand{\bx}{\strut\hfill$\blacksquare$\medbreak}

\newcommand{\ca}{\mathrm{C}^*}

\newcommand{\ol}{\overline}

\newenvironment{spmatrix}{\left(\begin{smallmatrix}}{\end{smallmatrix}\right)}

\newcommand{\bbC}{{\mathbb{C}}}

\newcommand{\bbZ}{{\mathbb{Z}}}
 \newcommand{\A}{{\mathcal{A}}}
 \newcommand{\B}{{\mathcal{B}}}
 \newcommand{\C}{{\mathcal{C}}}
 
 \newcommand{\E}{{\mathcal{E}}}

\renewcommand{\H}{{\mathcal{H}}}

 \newcommand{\M}{{\mathcal{M}}}
 
\renewcommand{\O}{{\mathcal{O}}}

 \newcommand{\R}{{\mathcal{R}}}

 \newcommand{\U}{{\mathcal{U}}}

\newcommand{\upchi}{{\raise.35ex\hbox{$\chi$}}}
\newcommand{\fE}{{\mathfrak{E}}}

\newcommand{\qand}{\quad\text{and}\quad}

\newcommand{\qfor}{\quad\text{for}\quad}

\newcommand{\qforal}{\quad\text{for all}\quad}

\newcommand{\Alg}{\operatorname{Alg}}
\newcommand{\diag}{\operatorname{diag}}

\newcommand{\ran}{\operatorname{Ran}}

\newcommand{\spn}{\operatorname{span}}

\newcommand{\Tr}{\operatorname{Tr}}

\newcommand{\fix}{\operatorname{Fix}}
\def\bra#1{\langle #1|}
\def\ket#1{|#1 \rangle}
\def\one{{\mathchoice{\rm 1\mskip-4mu l}{\rm 1\mskip-4mu l}{\rm 1\mskip-4.5mu l}{\rm
1\mskip-5mu l}}}

\newcommand{\fixed}{\fix (\Phi)}

\newcommand{\bofh}{\B(\H)}

\begin{document}

\title[Quantum Computing]%
{A Quantum Computing Primer for Operator Theorists}
%
\author[D.W.Kribs]{David~W.~Kribs}
\thanks{{\it key words and phrases. } quantum computation, quantum information, quantum algorithms, completely positive maps, quantum channels, quantum error correction, noiseless
subsystems.}
\address{Department of Mathematics and Statistics, University of Guelph,
Guelph, ON, CANADA  N1G 2W1. }
\address{Institute for Quantum Computing, University of
Waterloo, Waterloo, ON, CANADA N2L 3G1. }
\address{Perimeter Institute for Theoretical Physics, 35 King
St. North, Waterloo, ON, CANADA N2J 2W9.}

\email{dkribs@uoguelph.ca}
%
\begin{abstract}
This is an exposition of some of the aspects of quantum
computation and quantum information that have connections with
operator theory. After a brief introduction, we discuss quantum
algorithms. We outline basic properties of quantum channels, or
equivalently, completely positive trace preserving maps.  The main
theorems for quantum error detection and correction are presented
and we conclude with a description of a particular passive method
of quantum error correction.
\end{abstract}
\maketitle

\section{Introduction}\label{S:intro}

The study of the underlying mathematics for quantum computation
and quantum information is quickly becoming an interesting area of
research \cite{NC}. While these fields promise far reaching
applications \cite{Brooks,Johnson,Nielsen2002}, there are still
many theoretical and experimental issues that must be overcome,
and many  involve deep mathematical problems.  The main goal of
this paper is to provide a primer on some of the basic aspects of
quantum computing for researchers with interests in operator
theory or operator algebras. However, we note that the only
prerequisite for reading this article is a strong background in
linear algebra.

This work should not be regarded as an extensive introduction to
the subject. Indeed, the reader with knowledge of quantum
computing will surely have complaints about the selection of
material presented. Moreover, we do not consider here the
increasingly diverse fields of mathematics for which there are
connections with quantum computing (algebraic geometry, Fourier
analysis, group theory, number theory, operator algebras, etc).
Rather, our intention is to give a brief introduction and prove
some specific results with the hope that this paper will help
stimulate interest within the operator community.

The paper is organized as follows. The next section ($\S 2$)
contains a discussion of the basic notions, notation and
nomenclature used in quantum computing. In $\S 3$ we give a brief
introduction to quantum algorithms by describing some elementary
examples and presenting a simple algorithm (Deutsch's algorithm
\cite{Deutsch,DJ}) that demonstrates the power of quantum
computation. In $\S 4$ we outline the mathematical formalism for
the evolution of information inside quantum systems. This is
provided by quantum channels, which are represented by completely
positive trace preserving maps
\cite{Choi,Kraustext,Kraus,Paulsentext,Paulsentext2}. The
penultimate section ($\S 5$) includes a discussion of quantum
error correction methods. We present the fundamental theorems for
quantum error detection and correction in the `standard model'
\cite{Gottesman_IQEC,IntroQEC,Preskillplenum}. In the final
section ($\S 6$) we describe a specific method of quantum error
prevention \cite{DG,Lnoiseless,HKL,IntroQEC,KLV,LCW,ZR} called the
`noiseless subsystem via noise commutant' method. Finally, we have
included a large collection of references as an attempt to give
the interested reader an entrance point into the quantum
information literature.

\section{Quantum Computing Basics}\label{S:prelim}

Let $\H$ be a (complex) Hilbert space. We shall use the Dirac
notation for vectors and vector duals in $\H$: A typical vector in
$\H$ will be written as a `ket' $\ket{\psi}$, and the linear
functional on $\H$ determined by this vector is written as the
`bra' $\bra{\psi}$.  Notice that the products of a bra and a ket
yield the inner product, $\bra{\psi_1}\ket{\psi_2}$, and a rank
one operator, $\ket{\psi_2}\bra{\psi_1}$. In particular, given
$\ket{\psi}\in\H$, the rank one projection of $\H$ onto the
subspace $\{\lambda \ket{\psi}: \lambda\in\bbC \}$ is written
$\ket{\psi}\bra{\psi}$. Further let $\bofh$ be the collection of
operators which act on $\H$. We will use the physics convention
$U^\dagger$ for the adjoint of an operator $U$.

The study of operators on Hilbert space is central to the theory
of quantum mechanics. For instance, consider the following
formulation of the {\it postulates of quantum mechanics}
\cite{NC,vN}:

\vspace{0.05in}

$(i)$ To every closed quantum system there is an associated
Hilbert space $\H$. The state of the system at any given time is
described by a unit vector $\ket{\psi}$ in $\H$, or equivalently
by a rank one projection $\ket{\psi}\bra{\psi}$. When the state of
the system is not completely known, it is represented by a {\it
density operator} $\rho$ on $\H$, which is a positive operator
with trace equal to one. (This is the quantum analogue of a
probability distribution.)

$(ii)$ The notion of {\it evolution} in a closed quantum system is
described by unitary transformations. That is, there is a unitary
operator $U$ on the system Hilbert space such that the
corresponding evolution is described by the conjugation map $\rho
\mapsto U\rho U^\dagger$.

$(iii)$ A {\it measurement} of a quantum system on $\H$ is
described by a set of operators $M_k$, $1\leq k \leq r$, such that
\[
\sum_{k=1}^r M_k^\dagger M_k =\one.
\]
The measurement is {\it projective} if each of the $M_k$ is a
projection, and thus the $M_k$ have mutually orthogonal ranges. (A
`classical measurement' arises when all the projections are rank
one.) The index $k$ refers to the possible measurement outcomes in
an experiment. If the state of the system is $\ket{\psi}$ before
the experiment, then the probability that event $k$ occurs is
given by $p(k) \equiv \bra{\psi} M_k^\dagger M_k \ket{\psi}$.
Notice that $\{ p(k) \}$ determines a probability distribution.

$(iv)$ Given  Hilbert spaces $\H_1, \ldots, \H_m$ associated with
$m$ quantum systems, there is a {\it composite quantum system} on
the Hilbert space $\H_1\otimes \cdots \otimes \H_m$. In
particular, if the states of the individual systems are
$\ket{\psi_i}$, then the state of the composite system is given by
$\ket{\psi_1}\otimes \cdots \otimes \ket{\psi_m}$.

\vspace{0.05in}

The Hilbert spaces of interest in quantum information theory are
of dimension $N = d^n$ for some positive integers $n\geq 1$ and
$d\geq 2$. (It is generally thought that extensions to infinite
dimensional space will be necessary in the future, but the current
focus is mainly on finite dimensional aspects.) For brevity we
shall focus on the $d=2$ cases. Thus we let $\H_N$ be the Hilbert
space of dimension $2^n$ given by the $n$-fold tensor product
$\H_N = \bbC^2 \otimes \ldots \otimes \bbC^2 = (\bbC^2)^{\otimes
n}$. We will drop reference to $N$ when convenient. Let $\{
\ket{0}, \ket{1}\}$ be a fixed orthonormal basis for 2-dimensional
Hilbert space $\H_2 = \bbC^2$. These vectors will correspond to
the classical base states in a given two level quantum system;
such as the ground and excited states of an electron in an atom,
`spin-up' and `spin-down' of an electron, two polarizations of a
photon of light, etc. We shall make use of the abbreviated form
from quantum mechanics for the associated standard orthonormal
basis for $\H_{2^n} = (\bbC^2)^{\otimes n} \simeq \bbC^{2^n}$. For
instance, the basis for $\H_4$ is given by $ \{ \ket{ij} :
i,j\in\bbZ_2 \}$, where $\ket{ij}$ is the vector tensor product
$\ket{ij} \equiv \ket{i}\ket{j} \equiv \ket{i}\otimes\ket{j}$.

A quantum bit of information, or a `qubit', is given by a unit
vector $\ket{\psi} = a \ket{0} + b \ket{1}$ in $\H_2$. The cases
$a=0$ or $b=0$ correspond to the classical states, and otherwise
$\ket{\psi}$ is said to be in a {\it superposition} of the states
$\ket{0}$ and $\ket{1}$. A `qudit' is a unit vector in $\bbC^{d}$.
A vector state $\ket{\psi}$ in $\H_N$ is said to be {\it
entangled} if it cannot be written as a tensor product of states
from its component systems, so that $\ket{\psi}$ does not
decompose as $\ket{\psi} = \ket{\phi_1}\ket{\phi_2}$ for some
vectors $\ket{\phi_i}$, $i=1,2$. As an example, consider
$\ket{\psi} = \frac{\ket{00} + \ket{11}}{\sqrt{2}}$ in $\H_4$,
this is a vector from the so-called `EPR pairs' or `Bell states'.
Roughly speaking, the notion of {\it decoherence} in a quantum
system corresponds to the vanishing of off-diagonal entries in
matrices associated with the system as it evolves.

A number of specific unitary matrices arise in the discussions
below. The {\it Pauli matrices} are given by
\[
X = \left(\begin{matrix}
0&1 \\
1&0
\end{matrix}\right), \quad
Y =  \left(\begin{matrix}
0&-i \\
i&0
\end{matrix}\right), \quad
Z = \left(\begin{matrix}
1&0 \\
0&-1
\end{matrix}\right).
\]
Let $\one_2$ be the $2\times 2$ identity matrix. We regard these
as matrix representations for operators acting on the given basis
for $\H_2$. In the $n$-qubit case (so $N = 2^n$) we may consider
all `single qubit unitary gates' determined by the Pauli matrices.
This is the set of all unitaries $\{ X_k, Y_k, Z_k : 1 \leq k \leq
n\}$ where $X_1 = X \otimes \one_2^{\otimes (n-1)}$, $X_2 = \one_2
\otimes X  \otimes \one_2^{\otimes (n-2)}$, etc. Further let
$U_{CN}$ be the `controlled-NOT gate' on $\H_4$. This is the
unitary which acts on $\H_4$ by $ U_{CN} : \ket{i}\ket{j} \mapsto
\ket{i} \ket{(i+j)\, {\rm mod}\,2}. $ The CNOT gate has natural
extensions $\{ U_{CN}^{(k,l)} : 1\leq k\neq l \leq n\}$ to unitary
gates on $\H_N$, where the $k$th and $l$th tensor slots act,
respectively, as the control and target qubits. For instance, with
this notation $U_{CN} = U_{CN}^{(1,2)}$. Note that
$U_{CN}^{(k,l)}$ only acts on the $k$th and $l$th qubits in $\H_N
= (\bbC^2)^{\otimes n}$. The set $\{ X_k, Y_k, Z_k,U_{CN}^{(k,l)}:
1 \leq k\neq l \leq n \}$ forms a set of universal quantum gates
for $\H_N$, meaning that this set generates the set of $N\times N$
unitary matrices $\U(N)$ as a group (up to complex phases). We
shall also make use of the {\it Hadamard gate}
\[
H = \frac{1}{\sqrt{2}} \left(\begin{matrix} 1 & 1 \\ 1 & -1
\end{matrix}\right),
\]
and the spin-$\frac{1}{2}$ Pauli matrices $\sigma_k = 1/2 K$ for
$k=x,y,z$.


\section{Quantum Algorithms}\label{S:qalgs}

Simply put, a {\it quantum algorithm} consists of an ensemble of
initial states $\rho$ which evolves under a unitary matrix $U$ to
a final density matrix $U\rho U^\dagger$. The famous factoring
algorithm of Shor \cite{Shoralg2,Shoralg1} and search algorithm of
Grover \cite{Grover} have received tremendous attention over the
last decade. As a proper treatment of these algorithms is beyond
the scope of this article, in this note we shall present the
Deutsch algorithm \cite{Deutsch} (and its generalization, the
Deutsch-Josza algorithm \cite{DJ}) since it is easily accessible
and gives a good illustration of the power of quantum computation.
In doing so, we shall give a description of two fundamental
classes of quantum algorithms: the simulation of a classical
function on a quantum computer and the algorithm for quantum
parallel computation.

Before continuing, let us illustrate a simple example of how an
operation such as {\it addition} may be performed with a quantum
algorithm. This will also allow us to establish some notation for
the discussion which follows. We shall identify the standard basis
vectors $\ket{i_1\cdots i_n}$ for $\H_{N}$ with the integers
$\{0,1,\ldots , 2^n -1\}$ via binary expansions. Then we may
define a unitary $U$ on $\H_N \otimes \H_N$ by $U\,\ket{x}\ket{y}
= \ket{x}\ket{x\oplus y}$ where the addition $x\oplus y$ is modulo
$N$. The corresponding quantum algorithm implements addition
modulo $N$. (Note that the CNOT gate is obtained in the case
$N=2$.)

Towards the Deutsch algorithm, let $H_n = H^{\otimes n}$ be the
$n$-fold tensor product of the Hadamard gate acting on $\H_N$.
Observe that $H_n \ket{0}^{\otimes n}$ is the uniform
superposition
\[
H_n \ket{0}^{\otimes n} =\frac{1}{\sqrt{2^n}} \sum_{x=0}^{2^n -1}
\ket{x}.
\]

Fix positive integers $k,m\geq 1$. Let $f: \bbZ_2^m \rightarrow
\bbZ_2^k$ be a function. Consider the space $\H_{m,k} = \H_{2^m}
\otimes \H_{2^k}$. Then $\H_{m,k}$ has basis $\ket{x}\ket{y} =
\ket{x}\otimes \ket{y}$ where $x\in\bbZ_2^m$ and $y\in\bbZ_2^k$.
Define a unitary map $U_f : \H_{m,k} \rightarrow \H_{m,k}$ by
\[
U_f : \ket{x}\ket{y} \mapsto \ket{x} \ket{y\oplus f(x)}.
\]
For a given $x\in\bbZ_2^m$, notice that the action of $U_f$ is to
permute the basis vectors $\{ \ket{x}\ket{y} : y\in\bbZ_2^k \}$.

\begin{note}
Observe that $U_f \ket{x}\ket{0} = \ket{x}\ket{f(x)}$ for all $x$.
Thus $U_f$ {\it simulates} $f$ on a quantum computer, and in this
sense any classical function can be performed on a quantum
computer.
\end{note}

Next we compute
\begin{eqnarray*}
U_f \big( (H_m \ket{0}^{\otimes m} )\otimes \ket{0}^{\otimes k}
\big) &=& U_f \Big( \frac{1}{\sqrt{2^m}} \sum_{x=0}^{2^m -1}
\ket{x} \otimes \ket{0} \Big) \\
&=&   \frac{1}{\sqrt{2^m}} \sum_{x=0}^{2^m -1}
U_f (\ket{x}\otimes \ket{0}) \\
&=&  \frac{1}{\sqrt{2^m}} \sum_{x=0}^{2^m -1} \ket{x} \otimes
\ket{f(x)}
\end{eqnarray*}
Hence an application of $H_m \otimes \one_{2^k}$ followed by $U_f$
yields a simultaneous parallel computation of $f$ on {\it all}
possible values of $x$. The corresponding `circuit-gate' diagram
for {\it quantum parallelism} for $f$ is given below. In such a
diagram, the `circuits' correspond to states from the component
systems (in this case $\ket{0}^{\otimes m}$ and $\ket{0}^{\otimes
k}$ from the systems $\H_{2^m}$ and $\H_{2^k}$ respectively). The
`gates', drawn as boxes, indicate unitary operators applied as the
system evolves with a left-to-right convention in the diagram (so
here $H_m\otimes\one_{2^k}$ is applied first, then $U_f$ is
applied).

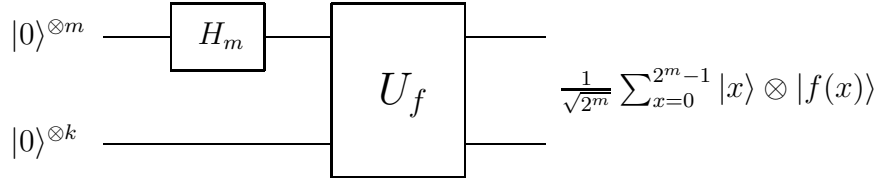
\setlength{\unitlength}{.007in}

\begin{picture}(100,200)(10,20)
\put(20,50){$\ket{0}^{\otimes k}$}

\put(20,130){$\ket{0}^{\otimes m}$}

\put(430,90){$\frac{1}{\sqrt{2^m}} \sum_{x=0}^{2^m -1}
\ket{x}\otimes\ket{f(x)}$}

\put(90,55){\line(1,0){170}}

\put(360,55){\line(1,0){60}}

\put(90,135){\line(1,0){50}}

\put(210,135){\line(1,0){50}}

\put(360,135){\line(1,0){60}}


\put(260,30){\line(1,0){100}}

\put(360,30){\line(0,1){130}}

\put(260,30){\line(0,1){130}}

\put(260,160){\line(1,0){100}}

\put(295,85){{\Large $U_f$}}


\put(140,110){\line(1,0){70}}

\put(140,110){\line(0,1){50}}

\put(210,110){\line(0,1){50}}

\put(140,160){\line(1,0){70}}

\put(160,130){$H_m$}

\end{picture}

Let $f:\bbZ_2 \rightarrow \bbZ_2$ be a function. Then call $f$
{\it constant} if $f(0) = f(1)$ and {\it balanced} if $f(0)\neq
f(1)$. The problem addressed by {\it Deutsch's algorithm} is the
following: Given $f:\bbZ_2 \rightarrow \bbZ_2$, determine if it is
constant or balanced. Classically an answer to this problem
requires two evaluations of $f$. On the other hand, Deutsch's
algorithm shows how to do this with one quantum operation. The
circuit-gate diagram for the algorithm is given below.

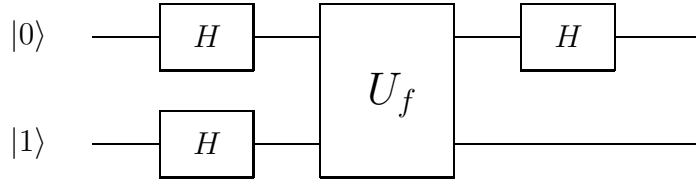
\begin{picture}(100,200)(-40,20)
\put(28,50){$\ket{1}$}

\put(28,130){$\ket{0}$}


\put(90,55){\line(1,0){50}}

\put(210,55){\line(1,0){50}}

\put(360,55){\line(1,0){180}}

\put(90,135){\line(1,0){50}}

\put(210,135){\line(1,0){50}}

\put(360,135){\line(1,0){50}}

\put(480,135){\line(1,0){60}}


\put(260,30){\line(1,0){100}}

\put(360,30){\line(0,1){130}}

\put(260,30){\line(0,1){130}}

\put(260,160){\line(1,0){100}}

\put(295,85){{\Large $U_f$}}


\put(140,110){\line(1,0){70}}

\put(140,110){\line(0,1){50}}

\put(210,110){\line(0,1){50}}

\put(140,160){\line(1,0){70}}

\put(165,128){$H$}


\put(140,30){\line(1,0){70}}

\put(140,30){\line(0,1){50}}

\put(210,30){\line(0,1){50}}

\put(140,80){\line(1,0){70}}

\put(165,48){$H$}


\put(410,110){\line(1,0){70}}

\put(410,110){\line(0,1){50}}

\put(480,110){\line(0,1){50}}

\put(410,160){\line(1,0){70}}

\put(435,128){$H$}


\end{picture}

\vspace{0.2in}

The initial state is given by $\ket{01} = \ket{0}\otimes\ket{1}$.
The first stage of the algorithm evolves this state to
\[
(H\otimes H) (\ket{0}\otimes \ket{1}) = (H\ket{0})\otimes
(H\ket{1}) = \ket{+}\otimes \ket{-},
\]
where $\ket{+} = \frac{\ket{0}+\ket{1}}{\sqrt{2}}$ and $\ket{-} =
\frac{\ket{0}-\ket{1}}{\sqrt{2}}$. Observe that the action of
$U_f$ as defined above yields
\[
U_f (\ket{x} \otimes \ket{-}) = (-1)^{f(x)} \ket{x} \otimes
\ket{-} \qfor x\in\bbZ_2.
\]
Thus, after the second stage of the algorithm the system has
evolved to
\begin{eqnarray*}
U_f (\ket{+}\otimes \ket{-}) &=& \frac{1}{\sqrt{2}} U_f \big(
(\ket{0} +
\ket{1}) \otimes \ket{-} \big) \\
&=& \Big( \frac{(-1)^{f(0)}\ket{0} +
(-1)^{f(1)}\ket{1}}{\sqrt{2}}\Big) \otimes \ket{-} \\
&=& \left\{ \begin{array}{cl} \pm \ket{+}\otimes\ket{-} & \mbox{if
$f$ is constant} \\
\pm \ket{-}\otimes\ket{-} & \mbox{if $f$ is balanced}
\end{array}\right.
\end{eqnarray*}
Finally, we apply a Hadamard gate $H$ to the first qubit; that is,
we apply $H\otimes \one_2$ to the full system. Thus the system
evolves to
\begin{eqnarray*}
\left\{\begin{array}{cl} (H\otimes \one_2)(\pm
\ket{+}\otimes\ket{-}) =
\pm \ket{0}\otimes\ket{-} & \mbox{if $f$ is constant} \\
(H\otimes \one_2)(\pm \ket{-}\otimes\ket{-}) = \pm
\ket{1}\otimes\ket{-} & \mbox{if $f$ is balanced}
\end{array}\right.
\end{eqnarray*}
In particular, if we measure the first qubit we get
\begin{eqnarray*}
\left\{\begin{array}{cl} \pm \ket{0} & \mbox{if $f$ is constant} \\
\pm \ket{1} & \mbox{if $f$ is balanced}
\end{array}\right.
\end{eqnarray*}
Note that there is no uncertainty in the result: If we measure the
first qubit and obtain $\pm \ket{0}$ (respectively $\pm \ket{1}$),
then we know $f$ is constant (respectively balanced) with
probability 1.

\begin{rem}
The Deutsch-Josza generalization \cite{DJ} yields a more dramatic
result. Let $f: \bbZ_2^m \rightarrow \bbZ_2$ be a function. Define
$f$ to be constant if $f(x) = f(y)$ for $x,y\in\bbZ_2^m$, and
balanced if $|f^{-1}(0)| = 2^{m-1} = |f^{-1}(1)|$. Suppose we know
$f$ is either constant or balanced and that we wish to determine
which one $f$ is. Classically, this requires $2^{m-1}+1$
evaluations to know with certainty what type $f$ is. The
Deutsch-Josza algorithm shows how this may be accomplished with a
single operation on a quantum computer. The algorithm acts on
$\H_{2^m}\otimes \H_2$ with initial state $\ket{0}^{\otimes m}
\otimes \ket{1}$. The circuit-gate diagram may be obtained by
adjusting the Deutsch diagram with $\ket{0}^{\otimes m} \otimes
\ket{1}$ in place of $\ket{0}\otimes \ket{1}$, $m$ circuits at the
top instead of one, $H_m$ in place of $H$ in the top circuits
prior to $U_f$, and $H_m$ in place of the final $H$. Thus,
\[
\ket{0}^{\otimes m} \otimes \ket{1} \mapsto (H_m\otimes \one_2)
U_f (H_m \otimes H)( \ket{0}^{\otimes m} \otimes \ket{1}).
\]
\end{rem}

\begin{note}
As a starting point for more recent work on quantum algorithms we
mention \cite{BHMT,MZ,Watrous}. Also see \cite{BBDR} for an
extensive mathematical introduction to the study of quantum
algorithms.
\end{note}





\section{Quantum Channels}\label{S:channels}

While evolution in a closed quantum system is unitary (see
postulate $(ii)$), experimentally evolution occurs in an `open
quantum system'. In such a system evolution is described
mathematically by completely positive trace preserving maps
\cite{NC}. The physical motivation for this is discussed below.

A (linear) map $\E: \bofh \rightarrow \bofh$ where $\H$ is a
Hilbert space is {\it completely positive} if for all $k\geq 1$
the `ampliation map'
\[
\E^{(k)} : \M_k (\bofh) \rightarrow \M_k (\bofh)
\]
given by $\E^{(k)} ( (\rho_{ij})) = (\E(\rho_{ij}))$ is a positive
map. (We could also write $\E^{(k)} = \one_k \otimes \E$.) This is
a rather strong condition; for instance, the transpose map is the
standard example of a positive map which is not completely
positive. The study of completely positive maps has been an active
area of research in both the quantum physics and operator theory
communities for at least thirty years. (In fact, it seems that
many general results on completely positive maps have been
obtained in the two fields without knowledge of the other.) See
the texts of Kraus \cite{Kraustext} and Paulsen
\cite{Paulsentext,Paulsentext2} for good treatments of the subject
from the two perspectives.

Thus, in the  general setting quantum information evolves through
an open quantum system via completely positive trace preserving
maps.

\begin{defn}
A {\it quantum channel} $\E : \bofh \rightarrow \bofh$ is a map
which is completely positive and trace preserving.
\end{defn}

Trace preservation for a channel is equivalent to requiring the
preservation of probabilities as  states evolve through a quantum
system. A channel is positive since density operators must evolve
to density operators, and it is completely positive because this
property must be preserved when the initial system is tensored
with other systems (as part of a composite system). Physically, an
open system can be regarded as lying inside a larger closed
quantum system where all evolution occurs in a unitary manner.
Thus, evolution in the open system can be regarded as a
`compression' of the unitary evolution on the larger closed
system. The mathematical formalism for this physical description
is provided by Stinespring's dilation theorem \cite{Stine}.

The following fundamental result for completely positive maps was
proved independently by Choi \cite{Choi} and Kraus \cite{Kraus}.
We present Choi's short  operator proof.

\begin{thm}\label{opsum}
Let $\E: \B(\H_N) \rightarrow \B(\H_N)$ be a completely positive
map. Then there are operators $E_k\in\B(\H_N)$, $1\leq k \leq
N^2$, such that
\begin{eqnarray}\label{opsumrepn}
\E (\rho) = \sum_{k=1}^{N^2} E_k \rho E_k^\dagger \qforal \rho \in
\B(\H_N).
\end{eqnarray}
\end{thm}

\Prf Let $e_{ij} =\ket{i}\bra{j}$ be the matrix units associated
with the standard basis for $\H_N$. Let $R = \E^{(N)} ((e_{ij}))$.
This matrix is positive by the $N$-positivity of $\E$. (In fact,
Choi proved that the positivity of $R$ characterizes complete
positivity of $\E$.) Consider a decomposition $R =
\sum_{k=1}^{N^2} \ket{a_k}\bra{a_k}$, where
$\ket{a_k}\in\bbC^{N^2}$ are (appropriately normalized)
eigenvectors for $R$. Let $\{P_i : 1\leq i \leq N\}$ be the family
of rank $N$ projections on $\bbC^{N^2}$ which have mutually
orthogonal ranges and satisfy $P_i R P_j = \E(e_{ij})$. Then
$\ket{a_k} = \sum_{i=1}^N P_i \ket{a_k}$. Define operators $E_k
:\bbC^N \rightarrow \bbC^N$ by $E_k \ket{i} \equiv P_i \ket{a_k}$,
so that
\[
R = \sum_k \sum_{i,j} P_i \ket{a_k}\bra{a_k} P_j = \sum_{i,j} P_i
\Big( \sum_k E_k \ket{i}\bra{j} E_k^\dagger \Big) P_j.
\]
Hence,
\[
\E(e_{ij}) = \E(\ket{i}\bra{j}) = P_i R P_j = \sum_{k=1}^{N^2} E_k
\ket{i}\bra{j} E_k^\dagger,
\]
and equation (\ref{opsumrepn}) holds by the linearity of $\E$. \bx

The decomposition (\ref{opsumrepn}) is referred to as the {\it
operator-sum representation} of $\E$ in quantum information
theory. The operators $E_k$ are called the {\it noise operators}
or {\it errors} of the channel. Note that trace preservation of a
channel is equivalent to its noise operators $E_k$ satisfying
\[
\sum_k E_k^\dagger E_k = \one.
\]

\begin{rem}
Choi's proof of this theorem provided the impetus for a recent
application of Leung \cite{Leungtomog} to `quantum process
tomography'. This is a procedure by which an unknown quantum
channel can be fully recovered from experimental data. To
accomplish this for a given channel, (\ref{opsumrepn}) shows that
it is enough to recover the noise operators $E_i$. As Leung
observes, Choi's proof shows the only experimental data required
to fully reconstruct the $E_i$ are the output states
$\E(\ket{i}\ket{j})$. This fact is the core of the tomography
procedure outlined in \cite{Leungtomog}.
\end{rem}

The following result, which we state without proof, shows
precisely how different sets of noise operators for the same
channel are related.

\begin{prop}\label{linearcombos}
Let $\{E_1, \ldots, E_r\}$ and $\{E_1^\prime, \ldots,
E_r^\prime\}$ be the noise operators for channels $\E$ and
$\E^\prime$ respectively. Then $\E = \E^\prime$ if and only if
there is an $r\times r$ scalar unitary matrix $U = (u_{ij})$ such
that $E = U E^\prime$ where $E^t = [E_1 \cdots E_r]$ and
$(E^\prime)^t = [E_1^\prime \cdots E_r^\prime]$. In other words,
\[
E_i = \sum_{j=1}^r u_{ij} E_j^\prime \qfor 1\leq i \leq r.
\]
\end{prop}

Let $\A=\Alg\{E_i,E_i^\dagger\}$ be the algebra generated by the
$E_i$ and $E_i^\dagger$. This is the set of  polynomials in the
$E_i$ and $E_i^\dagger$. An application of the Cayley-Hamilton
theorem on minimal polynomials from linear algebra shows that all
such polynomials may be written as polynomials with degree below
some global bound. In quantum computing, $\A$ is called the {\it
interaction algebra} for the channel. It is $\dagger$-closed by
definition, hence it is a finite dimensional $\ca$-algebra
\cite{Arvinvite,byeg,Tak}. Observe that, as a direct consequence
of the previous result, $\A$ can be seen to be a relic of the
channel; in other words, it is independent of the choice of noise
operators which satisfy (\ref{opsumrepn}) for the channel. This is
most succinctly seen in the case of unital channels, see
$\S$~\ref{S:noiseless} for details.

\begin{note}
We mention that an interesting and highly active area of current
research in quantum information theory revolves around the study
of quantum channel capacities. Specifically, there are a number of
deep mathematical problems which are concerned with computing the
capacity of a quantum channel to carry classical or quantum
information. The following references give a starting point into
the literature
\cite{Bennett,DS,HW,Holevo_coding,Holevo_capthm,HSR,HHHLT,King_depolar,Lloyd,SW,Shor_equiv,Winter}.
\end{note}

Note that the quantum measurement process (postulate $(iii)$)
naturally determines a quantum channel. Let us consider more
specific examples. (The final section includes further examples.)

\begin{egs}
$(i)$ Let $0 < p < 1$ and define operators on $\H_2$ by
\[
E_1 =  (\sqrt{1-p}) \one_2, \quad\quad E_2 = (\sqrt{p}) X,
\]
with respect to $\{\ket{0},\ket{1}\}$. The {\it bit flip} channel
$\E(\rho) = E_1 \rho E_1^\dagger + E_2 \rho E_2^\dagger$ flips the
state $\ket{0}$ to $\ket{1}$ and vice versa with probability $p$.
For instance, observe that the projection $\ket{0}\bra{0}$ evolves
to the superposition $\E(\ket{0}\bra{0}) = (1-p) \ket{0}\bra{0} +
p \ket{1}\bra{1}.$

$(ii)$ Let $E_1 = \frac{\one_2}{2}$ and define a channel by
$\E(\rho) = E_1\rho E_1 + \sigma_x\rho \sigma_x + \sigma_y\rho
\sigma_y +\sigma_z\rho \sigma_z$. Observe that $\E(\rho) =
\frac{\one_2}{2}$ for every density operator $\rho\in\M_2$. This
property, distinct density matrices evolving in a channel to the
same density matrix, has recently been exploited as part of a
scheme for quantum cryptography \cite{BRS1}.

$(iii)$ Let $0< r < 1$ and define operators on $\H_2$ by
\[
E_1 =  \left( \begin{matrix}
1 & 0 \\
0 & \sqrt{1-r}
\end{matrix}\right), \quad\quad E_2 =  \left( \begin{matrix}
0 & \sqrt{r} \\
0 & 0
\end{matrix}\right).
\]
These noise operators define an {\it amplitude damping} channel
$\E(\rho) = E_1 \rho E_1^\dagger + E_2 \rho E_2^\dagger$. Such
channels characterize energy dissipation within a quantum system.

$(iv)$ Let $U_1, \ldots, U_d$ be unitaries which act on a common
Hilbert space and let $r_1,\ldots ,r_d$ be positive scalars such
that $\sum_i r_i = 1$. Then we may define a channel by $\E(\rho) =
\sum_{i=1}^d r_i U_i \rho U_i^\dagger$. These are the prototypical
examples of `unital' ($\E(\one)=\one$) channels (see
$\S$~\ref{S:noiseless}). In fact, every unital channel on $\M_2$
may be written as a convex sum of unitaries in this way. This is
not the case in higher dimensions however.

$(v)$ The class of {\it entanglement breaking channels} was
introduced in \cite{Holevo_coding} and studied in \cite{HSR}.
These are quantum channels which can be written in the form
\[
\E(\rho) = \sum_k \ket{\psi_k}\bra{\psi_k} \bra{\phi_k}\rho
\ket{\phi_k},
\]
for some vectors $\ket{\psi_k}$ and $\ket{\phi_k}$. With this
representation, trace preservation is equivalent to $\sum_k
\ket{\phi_k}\bra{\phi_k} = \one$, as $\Tr (\rho
\ket{\phi_k}\bra{\phi_k}) = \bra{\phi_k}\rho \ket{\phi_k}$. Such
channels derive their name from the fact that for $d\geq 1$, a
density operator $ \E^{(d)} (\Gamma)$ is never entangled, even if
$\Gamma$ was initially entangled.
\end{egs}

\section{Quantum Error Correction}\label{S:qec}

In this section we present the central aspects of the `standard
model' for quantum error detection and correction. For an
extensive introduction to the subject we point the reader to the
articles \cite{Gottesman_IQEC,IntroQEC,Preskillplenum}. The
general error correction problem in quantum computing is much more
delicate when compared to  error correction  in classical
computing. The possible errors that can occur include all possible
unitary matrices, whereas in classical computing the only errors
are bit flips. Nonetheless, methods have been (and are being)
developed for quantum error correction.

\subsection{Error Detection}\label{sS:prevent}

Let $\H$ be the Hilbert space for a given quantum system. Then a
{\it quantum code} $\C$ on $\H$ is a subspace of $\H$. Let $P_\C$
be the projection of $\H$ onto $\C$. Then $P_\C$ and $P_\C^\perp$
describe a measurement of the system which can be used to
determine if a given state $\ket{\psi}\in\H$ belongs to the code.
The basic idea of an error-detection scheme in this setting is to
first prepare an initial state in $\C$, for brevity let us
restrict ourselves to unit vectors $\ket{\psi}$ in $\C$. The state
$\ket{\psi}\bra{\psi}$ is then transmitted through the quantum
channel $\E$ of interest, evolving to $\E(\ket{\psi}\bra{\psi})$.
Finally, the measurement $P_\C,P_\C^\perp$ is performed on this
final state. This motivates the following definition.

\begin{defn}
Let $\C$ be a quantum code on $\H$ and let $E$ be an error (noise)
operator associated with a given quantum channel on $\H$. Then
$\C$ {\it detects the error $E$} if the states which are accepted
after $E$ acts are unchanged, up to a scaling factor. In other
words, there is a scalar $\lambda_E$ such that
\[
P_\C E \ket{\psi} = \lambda_E \ket{\psi} \qforal \ket{\psi}\in \C.
\]
\end{defn}

Observe that the set of error operators $E$ which are detectable
for a fixed quantum code $\C$ form a subspace of operators, or a
so-called operator space \cite{Paulsentext}. There are a number of
equivalent conditions for detectable errors.

\begin{thm}\label{qed}
Let $\C$ be a quantum code and let $E$ be an error operator
associated with a given quantum channel. Then the following
conditions are equivalent:
\begin{itemize}
\item[$(i)$] $E$ is detectable by $\C$, with scaling factor
$\lambda_E$. \item[$(ii)$] $P_\C E P_\C = \lambda_E P_\C$.
\item[$(iii)$] $\bra{\psi_1} E \ket{\psi_2} = \lambda_E
\bra{\psi_1} \ket{\psi_2}$ for all $\ket{\psi_i}\in \C$, $i=1,2$.
\item[$(iv)$] For every pair of vectors $\ket{\psi_1}$ and
$\ket{\psi_2}$ in $\C$ which are orthogonal, the vectors $E
\ket{\psi_1}$ and $\ket{\psi_2}$ are orthogonal.
\end{itemize}
\end{thm}

\Prf We shall prove the implication $(iv)\Rightarrow (iii)$. The
other directions are either trivial or easy to see. We may clearly
assume that $\dim \C \geq 2$. Let $\{ \ket{\psi_1}, \ket{\psi_2},
\ldots \}$ be an orthonormal basis for $\C$. We first claim that
$(iv)$ implies
\begin{eqnarray}\label{lambdae}
\bra{\psi_i} E \ket{\psi_i} = \bra{\psi_j} E \ket{\psi_j} \qforal
i,j.
\end{eqnarray}
Indeed, to see this, fix $i,j$ and define
\[
\ket{+} = \ket{\psi_i} + \ket{\psi_j} \qand \ket{-} = \ket{\psi_i}
- \ket{\psi_j}.
\]
Then by $(iv)$ we have
\[
0 = \bra{+} E \ket{-} = \bra{\psi_i}E \ket{\psi_i} - \bra{\psi_j}E
\ket{\psi_j}.
\]
Thus define $\lambda_E = \bra{\psi_i}E \ket{\psi_i}$, and note
that this is independent of $i$. Now let $\ket{\psi} = \alpha_1
\ket{\psi_1} + \alpha_2 \ket{\psi_2}+ \ldots$ and $\ket{\phi} =
\beta_1 \ket{\psi_1} + \beta_2 \ket{\psi_2}+ \ldots$ be vectors in
$\C$. Then by $(iv)$ and (\ref{lambdae}) we have
\[
\bra{\psi} E \ket{\phi} = \sum_{i,j} \alpha_i \ol{\beta_j}
\bra{\psi_i} E \ket{\psi_j} = \sum_i \alpha_i \ol{\beta_i}
\bra{\psi_i} E \ket{\psi_i} = \lambda_E \bra{\psi} \ket{\phi},
\]
and this completes the proof. \bx

Let us describe a matrix perspective for detectable errors. This
can be seen through a simple example.

\begin{eg}\label{3rep}
Let $\C$ be the quantum code given by
\[
\C = \spn \{ \ket{000}, \ket{111}\}\subseteq \H_8.
\]
The (unnormalized) error operators for the $n$-qubit depolarizing
channel \cite{NC,King_depolar} are $\fE = \{\one_{2^n},
Z_1,\ldots, Z_n \}$. In the 3-qubit case consider the error
operator $E \equiv Z_1$. Observe that $E\ket{000} = \ket{000}$
whereas $E\ket{111} = - \ket{111}$. Thus if $E$ was detectable by
$\C$ we would have,
\[
\lambda_E = \bra{000} E \ket{000} = 1 \qand \lambda_E = \bra{111}
E \ket{111} = -1.
\]

From this contradiction it follows  that none of the depolarizing
errors $Z_k$ are detectable by the code $\C$. In fact, from
Theorem~\ref{qed} the detectable errors for $\C$ may be realized
in matrix form as the operator space $\big\{\begin{spmatrix}
\lambda \one_2 & \ast \\ \ast & \ast
\end{spmatrix} : \lambda \in \bbC \big\},$ with respect to an
ordered orthonormal basis for $\H_8$ of the form $\{ \ket{000},
\ket{111}, \ldots\}$, and where the $\ast$ entries indicate that
any choice is admissible.
\end{eg}

More generally, the conditions of Theorem~\ref{qed} show that the
detectable errors for a given quantum code $\C$ form the subspace
of operators
\[
\left\{ \left(
\begin{matrix}
\lambda \one & \ast \\
\ast & \ast
\end{matrix}\right) : \lambda \in \bbC \right\},
\]
where the matrix form is given with respect to the spatial
decomposition $\H = P_\C\H \oplus P_\C^\perp \H$.

\subsection{Error Correction}\label{S:correct}

Let $\fE = \{ E_i \}$ be a set of errors that act as noise
operators for a given quantum channel. If $\C$ is a quantum code,
then the basic error-correction problem for $\C$ is to determine
when there is a decoding procedure for $\C$ such that all the
errors in $\fE$ are corrected. The simplest possible case occurs
when every error $E_i$ is the multiple of a unitary operator and
the subspaces $E_i\C$ are mutually orthogonal. The obvious
decoding procedure in this situation is to first make a projective
measurement to determine which of the subspaces $E_i \C$ a given
state $\ket{\psi}\in\C$ has evolved to, then apply the inverse of
the error operator $E_i$. This particular case motivates the
following general definition in the standard model for quantum
error correction.

\begin{defn}
Let $\E$ be a quantum channel and let $\C$ be a quantum code with
projection $P_\C$. Then $\C$ is {\it correctable for $\E$} if
there is a quantum channel $\R$ such that
\[
\R \circ \E (\rho) = \rho
\]
for all $\rho$ supported on $\C$; that is, all $\rho$ with $\rho =
P_\C \rho P_\C$.
\end{defn}

There are a number of useful characterizations of correctable
codes. Note that in the proof of $(iii)\Rightarrow (i)$ below, the
error correction operation $\R$ is explicitly constructed and
hence this gives a constructive approach for decoding within this
error correction model.

\begin{thm}\label{correctthm}
Let $\E$ be a quantum channel  with errors $\fE = \{ E_i \}$ and
let $\C$ be a quantum code with projection $P_\C$. Then the
following conditions are equivalent:
\begin{itemize}
\item[$(i)$] $\C$ is correctable for $\E$. \item[$(ii)$] The
operators in the set $\fE^\dagger \fE = \{ E_1^\dagger E_2 :
E_i\in\fE \}$ are detectable by $\C$. \item[$(iii)$] There are
scalars $\Lambda = (\lambda_{ij})$ such that
\begin{eqnarray}\label{correcteqn}
P_\C E_i^\dagger E_j P_\C = \lambda_{ij} P_\C \qforal i,j.
\end{eqnarray}
\item[$(iv)$] There is a linear transformation $E_i \mapsto
E_i^\prime$ on $\fE$ such that the new error operators
$E_i^\prime$ satisfy the following properties:
\begin{itemize}
\item[$(a)$] The subspaces $E_i^\prime \C$ are mutually
orthogonal, \item[$(b)$] The restriction of every $E_i^\prime$ to
$\C$ is proportional to a restriction to $\C$ of a unitary
operator.
\end{itemize}
\end{itemize}
\end{thm}

\Prf Observe that any matrix $\Lambda$ which satisfies
(\ref{correcteqn}) must be positive since this equation may be
written as a matrix product $A^\dagger A = (\lambda_{ij}P_\C)$
where $A=[E_1P_\C \,\, E_2P_\C \cdots]$ is a row matrix.
Conditions $(ii)$ and $(iii)$ are equivalent by definition. We
shall prove $(i)$ is necessary and sufficient for $(iii)$ and
leave the connection with $(iv)$ for the interested reader.

For $(i)\Rightarrow (iii)$, let $\C$ be a quantum code with code
projection $P_\C$. Suppose $\E$ is a quantum channel with errors
$\{E_i\}$ and that $\C$ is correctable for $\E$ via the
error-correction operation $\R$ with noise operators $\{R_j\}$.
Define a compressed channel by $\E_\C (\rho) \equiv \E(P_\C\rho
P_\C)$. Then by hypothesis $\R (\E_\C(\rho)) = \R (\E(P_\C\rho
P_\C)) = P_\C\rho P_\C$. In particular,
\[
\sum_{i,j} R_j E_i P_\C \rho P_\C E_i^\dagger R_j^\dagger = P_\C
\rho P_\C \qforal \rho.
\]
Thus by Proposition~\ref{linearcombos} there are scalars
$\alpha_{ki}$ such that
\[
R_k E_i P_\C = \alpha_{ki} P_\C \qforal i,k.
\]
Hence,
\begin{eqnarray}\label{Rk}
P_\C E_i^\dagger R_k^\dagger R_k E_j P_\C & = &\ol{\alpha}_{ki}
\alpha_{kj} P_\C \qforal i,j,k.
\end{eqnarray}
But $\R$ preserves traces, so that $\sum_k R_k^\dagger R_k = \one$, and thus when we sum
(\ref{Rk}) over $k$ we find
\[
P_\C E_i^\dagger E_j P_\C = \lambda_{ij} P_\C \qforal i,j,
\]
where $\lambda_{ij} = \sum_k  \ol{\alpha}_{ki} \alpha_{kj}$.

Conversely, let us assume that $\E$ is a channel with errors
$\{E_i \}$ and $\C$ is a code with projection $P_\C$ such that
(\ref{correcteqn}) holds for a positive scalar matrix $\Lambda =
\{\lambda_{ij}\}$. Let $U$ be a unitary operator such that $D =
U^\dagger \Lambda U = (d_{kl})$ is diagonal (so that $D=\diag
(d_{kk})$). Note that $\sum_k d_{kk} =1$ by trace preservation of
$\E$. By Proposition~\ref{linearcombos} the operators $F_k =
\sum_i u_{ik} E_i$ also implement $\E$. A simple computation shows
that
\[
P_\C F_k^\dagger F_l P_\C = d_{kl} P_\C \qforal k,l.
\]

The polar decomposition of $F_k P_\C$ gives
\[
F_k P_\C = U_k \sqrt{P_\C F_k^\dagger F_k P_\C} = \sqrt{d_k} U_k
P_\C
\]
for some unitary $U_k$. Define projections $P_k \equiv U_k P_\C
U_k^\dagger$. Then
\[
P_l P_k = P_l^\dagger P_k = \frac{U_l P_\C F_l^\dagger F_k P_\C
U_k^\dagger}{\sqrt{d_{ll}d_{kk}}} = 0 \qfor k\neq l,
\]
and hence the ranges of the $P_k$ are mutually orthogonal.

Without loss of generality assume that $\sum_k P_k = \one$
(otherwise just add the projection onto the orthogonal complement
and define $U_k = \one$). The candidate error-correction operation
is defined by
\[
\R (\rho) = \sum_k U_k^\dagger P_k \rho P_k U_k.
\]
Observe that for all $\rho$ with $\rho = P_\C \rho P_\C$,
\begin{eqnarray*}
U_k^\dagger P_k F_l \sqrt{\rho} = U_k^\dagger P_k^\dagger F_l \sqrt{\rho}
&=& \frac{U_k^\dagger U_kP_\C F_k^\dagger F_l P_\C \sqrt{\rho}}{\sqrt{d_{kk}}} \\
&=&\delta_{kl} \sqrt{d_{kk}} P_\C \rho   =\delta_{kl}
\sqrt{d_{kk}} \rho.
\end{eqnarray*}
Thus we have
\[
\R(\E(\rho)) = \sum_{k,l} U_k^\dagger P_k F_l \rho F_l^\dagger P_k U_k = \sum_{k,l}
\delta_{kl} d_{kk} \rho = \rho
\]
for all $\rho = P_\C \rho P_\C$, and we have proved that $(iii)
\Rightarrow (i)$. \bx



Let us discuss a pair of illustrative applications of
Theorem~\ref{correctthm}.

\begin{egs}
$(i)$ With $\C = \spn\{ \ket{000},\ket{111}\}$ inside $\H_8$ as in
Example~\ref{3rep}, it is easy to see that $\C$ satisfies
(\ref{correcteqn}) for the errors $\E =\{\one,X_1,X_2,X_3\}$. Let
$P_0$ be the projection onto $\C$, let $P_1$ be the projection
onto the subspace $X_1\C=\spn\{\ket{100},\ket{011}\}$ and
similarly define projections $P_2$ and $P_3$. Then the correction
operation given by the proof above is $\R = \{P_0,
X_1P_1,X_2P_2,X_3P_3\}$.


$(ii)$ Shor's 9-qubit code \cite{Shor95} is defined by two
orthonormal vectors in $\H_{2^9}$ given by
\begin{eqnarray*}
\ket{0_L} & = & \frac{(\ket{000} + \ket{111}) (\ket{000} +
\ket{111})(\ket{000} + \ket{111})  }{2\sqrt{2}} \\
\ket{1_L} & = & \frac{(\ket{000} - \ket{111}) (\ket{000} -
\ket{111})(\ket{000} - \ket{111})  }{2\sqrt{2}}.
\end{eqnarray*}
Let $\C = \spn\{ \ket{0_L}, \ket{1_L}\}$. Given $1\leq k \leq 9$,
the code $\C$ is correctable for the errors $\{ X_k, Y_k, Z_k\}$
as (\ref{correcteqn}) is satisfied for this triple. Let $P_{x,k}$
be the projection onto the subspace $X_k\C$ and similarly define
projections $P_{y,k}$ and $P_{z,k}$. Then the correction operation
given by the theorem is $\R = \{ X_k P_{x,k},\, Y_k P_{y,k},\, Z_k
P_{z,k} \}$.

Since the Pauli matrices, together with the identity operator
$\one_2$, form a linear basis for $\M_2$ that is closed under
multiplication up to scalar multiples, it follows that $\C$ is
correctable for any set of errors $\fE$ which act on one of the
nine possible qubits. (In fact the Shor code also corrects for
errors on multiple qubits \cite{NC}.)
\end{egs}

\begin{note}
The quantum error correction conditions of Theorems~\ref{qed} and
\ref{correctthm} were established independently by Bennett,
DiVincenzo, Smolin and Wooters \cite{BDSW96} and Knill and
Laflamme \cite{KL97}. As a collection of entry point references
for particular methods of quantum error correction we mention
\cite{CPZ96,Got97,Got96,Kit,Kit97b,KriLaf,Sho96,Ste96a,ZL96}.
\end{note}

\section{Noiseless Subsystems via the Noise Commutant}\label{S:noiseless}

In this section we describe a specific method of passive quantum
error correction, by which we mean no active intervention is
required after information is encoded. The basic idea in a
`noiseless subsystem method' of error correction, classical or
quantum, is to encode information on subsystems of the system of
interest in such a way that it remains immune to the effects of
the channel that the information is evolving through. Let $\E:
\bofh \rightarrow \bofh$ be a quantum channel with noise operators
$\{ E_i \}$ and interaction algebra $\A = \Alg \{E_i,
E_i^\dagger\}$. Recall that $\E$ is {\it unital} if $\E(\one) =
\sum_i E_i E_i^\dagger = \one$. The {\it noise commutant} for $\E$
is the $\dagger$-algebra
\begin{eqnarray*}
\A^\prime &=& \{ \rho \in\bofh : \rho A=A\rho \qfor A\in\A\} \\
&=& \{ \rho \in \bofh : [\rho,E_i]=0=[\rho,E_i^\dagger] \qfor
i=1,\ldots,n\}.
\end{eqnarray*}
The procedure in the {\it noiseless subsystem via noise commutant}
method of quantum error correction
\cite{DG,Lnoiseless,HKL,IntroQEC,KLV,LCW,Zan,ZR} is to use the
structure of $\A^\prime$ to produce noiseless subsystems (which
are also called `decoherence-free subspaces' in special cases, see
Remark~\ref{noiselessrem}).

The illustrations of this method that appear in the literature all
involve unital channels. The following discussion shows why this
is the case. Let
$
\fix (\E) = \{ \rho\in \bofh : \E(\rho) = \rho\}.
$
Observe that $\fix(\E)$ is a $\dagger$-closed subspace of $\bofh$.
Now consider a unital channel $\E$. Let $\rho $ belong to
$\A^\prime$. Then $\rho E_i = E_i\rho$ for all $i$ and hence
\[
\E(\rho) = \sum_i E_i \rho E_i^\dagger = \rho \E(\one) = \rho.
\]
Thus $\rho$ belongs to $\fix(\E)$. The converse inclusion was
proved independently in \cite{BS} and \cite{Kchannel} (see also
\cite{HKL} for another proof): The set $\fix(\E)$ coincides with
the noise commutant in the case of a unital channel. In fact, it
is not hard to show that the identity $\fix(\E) = \A^\prime$
characterizes unital channels. Thus, by building on the proofs
from \cite{BS,HKL,Kchannel} we may state the following.

\begin{thm}\label{unitalcondition}
Let $\E$ be a quantum channel. Then $\fix(\E) = \A^\prime$ if and
only if $\E$ is unital. Moreover, in this case $\A$ is equal to
the algebra $\A_0$ generated by the $E_i$. In other words, the
operators $E_i^\dagger$ belong to $\A_0$.
\end{thm}

\begin{note}
Notice how the algebra $\A = (\A)^{\prime\prime} =
\fix(\E)^\prime$ is a relic of the channel in the unital case (see
the discussion after Proposition~\ref{linearcombos}).
\end{note}

\begin{rem}\label{noiselessrem}
In the case of unital channels, Theorem~\ref{unitalcondition}
shows that the noiseless subsystem via noise commutant method can
be presented as follows: Let $\E$ be a unital quantum channel.
Then the interaction algebra $\A$ is generated by the $E_i$ and
$\A^\prime = \fix(\E)$ is a finite dimensional $\ca$-algebra. As
such, $\A^\prime$ is unitarily equivalent to a unique direct sum
of ampliated full matrix algebras, $\A^\prime \cong \oplus_k
(\one_{m_k}\otimes \M_{n_k})$.  Hence, a density operator $\rho$
that encodes the initial states of a quantum system will be immune
to the noise of the channel as it evolves through, $\E(\rho) =
\rho$, provided that $\rho$ is initially prepared on one of the
ampliated matrix blocks $\one_{m_k}\otimes \M_{n_k}$ inside the
noise commutant $\A^\prime = \fix(\E)$. The decoherence-free
subspace method may be regarded as the special case of this method
that occurs when matrix blocks $\one_{m_k}\otimes \M_{n_k}$ with
$m_k=1$ are utilized.
\end{rem}

In the case of a general quantum channel, however, the full
structure of $\A^\prime$ cannot be used for error correction. This
can be seen most dramatically in the following simple case.

\begin{prop}\label{annihilate}
Let $\E$ be a completely positive map such that $A_\E \equiv
\E(\one)$ is not invertible. Let $P_\E$ be the projection onto the
subspace $\H_\E \equiv \ran(A_\E)$. Suppose that $E_i = P_\E E_i
P_\E$ for all $i$. Then $\H_\E^\perp$ is non-zero and every
operator $\rho \in \B(\H_\E^\perp)$ belongs to $\A^\prime$ and
satisfies $\E(\rho) = 0$.
\end{prop}

\Prf The non-invertibility of $A_\E$ implies $\H_\E^\perp$ is
non-zero. Let $\rho \geq 0$ belong to $\B(\H_\E^\perp)$, by which
we mean $\rho$ acts on $\H$ with $\rho = P_\E^\perp \rho
P_\E^\perp$. Then $\rho$ trivially belongs to $\A^\prime$ since
$\rho E_i = 0 = E_i \rho$ for all $i$. Further, $\E(\rho) =
\E(P_\E^\perp \rho) = 0$. \bx

\begin{rem}
It would be interesting to know if the noise commutant can be used
to produce noiseless subsystems for classes of non-unital
channels. For instance, a natural generalization of unital
channels is the class of channels for which the identity evolves
to a multiple of a projection. Notice that if $\E(\one_\H) = m P$
for some projection $P$, then $m$ divides the dimension of $\H$ by
trace preservation of $\E$. A simple example of this phenomena is
given by the channel $\E$ with noise operators $A_i =
\ket{0}\bra{i}$ for $1\leq i \leq \dim\H\equiv d$. Trace
preservation of $\E$ may be readily verified, and in this case
\[
\E(\one_d) = \sum_{i=1}^d A_i A_i^\dagger = \sum_{i=1}^d
(\ket{0}\bra{i})(\ket{i}\bra{0}) = d \ket{0}\bra{0}.
\]
\end{rem}

\begin{note}
There are a number of other quantum error correction schemes that
have been investigated and some of these will be of interest to
operator researchers (See
\cite{AB,DG,Fi,KBLW,KriLaf,LCW,ZL,Zan,ZR}).
\end{note}

We finish by considering noiseless subsystems for some special
cases of unital channels.

\begin{egs}
$(i)$  Let $0 < p < 1$ and let $E_1, E_2$ be operators on $\H_2$
defined on the standard  basis by $E_1 = (\sqrt{1-p}) \one_2$ and
$E_2 = (\sqrt{p})Z$. Then $E_1$ and $E_2$ are the noise operators
for a unital channel $\E$ on $\M_2$ which is a variant on the bit
flip channel discussed earlier. The quantum operation
corresponding to this channel is equivalent to the {\it phase
flip} or {\it phase damping} operation on single qubits \cite{NC}.
It is so named because, for instance, $\E(\ket{+}\bra{+}) =
(1-p)\ket{+}\bra{+} + p \ket{-}\bra{-}$ and hence $\E$ flips the
phase of $\ket{+}\bra{+}$ and $\ket{-}\bra{-}$ with probability
$p$. It is easy to see in this case that
\begin{eqnarray*}
\fixed = \A^\prime &=& \{E_1,E_2\}^\prime \\
&=& \left\{   \left( \begin{matrix}
a & 0 \\
0 & b
\end{matrix}\right) : a,b\in\bbC \right\} \cong \bbC {\mathbf
1}\oplus \bbC {\mathbf 1},
\end{eqnarray*}
and hence this channel has no non-trivial noiseless subsystems.

$(ii)$ Let $0 < p < 1$ and let $E_1, E_2$ be operators on $\H_4$
defined on the standard basis by $E_1 = \sqrt{1-p} (\one_2 \otimes
\one_2)$ and $E_2 =  \sqrt{p}(Z\otimes Z)$. These noise operators
determine a unital channel $\E$ on $\M_4$ which can be regarded as
an ampliation of the phase flip channel. Compute
\begin{eqnarray*}
\fixed = \A^\prime &=& \{E_1,E_2\}^\prime  \\
 &=& \left\{   \left(
\begin{matrix}
a_{11} & 0 & 0 & a_{14} \\
0 & a_{22} & a_{23} & 0  \\
0 & a_{32} & a_{33} & 0  \\
a_{41} & 0 & 0 & a_{44}
\end{matrix}\right) : a_{ij} \in\bbC \right\}.
\end{eqnarray*}
Thus $\A^\prime$ is unitarily equivalent to the direct sum
$\A^\prime \cong \M_2 \oplus \M_2$ and there is a pair of
2-dimensional noiseless subsystems.

$(iii)$ More generally, let $E_1= U_1,\ldots ,E_d =U_d$ be
unitaries on a Hilbert space $\H$. If we are given scalars
$\lambda_i$ such that $\sum_i |\lambda_i|^2 = 1$, then the
operators $\lambda_i E_i$ are the noise operators for a unital
channel $\E$ with $\fix (\E) = \A^\prime = \{ U_1, \ldots,
U_d\}^\prime$. In this case the automatic self-adjointness of $\A$
can be seen directly through an application of the Cayley-Hamilton
theorem. (The algebra generated by any unitary $U$ can be seen to
include $U^\dagger$ via minimal polynomial considerations.)

$(iv)$ An important special case of such channels is the class of
`collective rotation channels'
\cite{BRS,LPS,Fi,Lnoiseless,HKL,KBLW,KLV,VKL,VFPKLC,ZL,Zan}. The
$n$-qubit example has noise operators given by weighted
exponentiations of the operators $J_k = \sum_{m=1}^n J_k^{(m)}$
for $k = x,y,z$, where $J_k^{(1)} = \sigma_k \otimes
(\one_2)^{\otimes (n-1)}$, etc. There is an abundance of noiseless
subsystems for these channels and they can be computed directly
\cite{HKLP} or through combinatorial techniques discussed below.

$(v)$ A generalization of the collective rotation class is
presented in \cite{JKK} and noiseless subsystems for this class of
`universal collective rotation (ucr) channels' are computed using
Young tableaux combinatorics. For each pair of positive integers
$d,n\geq 2$ there is a family of such channels with the case $d=2$
yielding the class from $(iv)$. It is proved in \cite{JKK} that
every ucr-channel possesses an abundance of noiseless subsystems.
In fact, it is shown that the noise commutant for every channel in
this class (for fixed $d$ and $n$) contains the algebra $\A_\pi =
\Alg\{\pi(\sigma):\sigma\in S_n\}$ where $S_n$ is the symmetric
group on $n$ letters and $\pi: S_n \rightarrow \U(d^n)$  is the
unitary representation of $S_n$ on $\H_{d^n}$ given by
\[
\pi(\sigma) \big( h_1\otimes \ldots \otimes h_n \big) =
h_{\sigma(1)} \otimes \ldots \otimes h_{\sigma(n)},
\]
for all $h_1,\ldots, h_n\in \H_d$ and $\sigma\in S_n$. In
particular, the Young tableaux machine can be used to explicitly
compute the structure of $\A_\pi$ and identify noiseless
subsystems for the corresponding ucr-channel. Further, a recent
paper of Bacon, et al \cite{BCH}, has shown that the change of
basis transformation from the standard basis to the Young tableaux
basis can be implemented efficiently with a quantum algorithm.
\end{egs}


{\noindent}{\it Acknowledgements.} I am grateful to the referees
for helpful comments. I would like to thank all the participants
in the Quantum Information Theory Learning Seminar organized by
the author and John Holbrook at the University of Guelph during
the fall of 2003. I am also grateful for enlightening
conversations with Daniel Gottesman, John Holbrook, Peter Kim,
Raymond Laflamme, Michele Mosca, Ashwin Nayak, Eric Poisson, David
Poulin, Martin Roetteler, Jean-Pierre Schoch, Robert Spekkens and
Paolo Zanardi. Support from NSERC, the University of Guelph, the
Institute for Quantum Computing and the Perimeter Institute is
kindly acknowledged.


\vspace{0.1in}



\end{document}